\documentclass[10pt,a4paper]{article}
\usepackage[utf8]{inputenc}

\usepackage{amsmath}
\usepackage{amsfonts}
\usepackage{amssymb}
\usepackage{url}
\usepackage{mathtools}
\usepackage{hyperref}

\title{The square of the 9-hypercube is 14-colorable\thanks{Work partially supported by the Emil Aaltonen Foundation}}
\author{Juho Lauri\thanks{Tampere University of Technology, Finland}}
\date{\today}

\begin{document}
\maketitle

\begin{abstract}
The $n$-hypercube, denoted by $Q_n$, has a vertex for each bit string of length $n$ with two vertices adjacent whenever their Hamming distance is one. The minimum number of colors needed to color $Q_n$ such that no two vertices at a distance at most $k$ receive the same color is denoted by $\chi_{\bar{k}}(n)$. Equivalently, $\chi_{\bar{k}}(n)$ denotes the minimum number of binary codes with minimum distance at least $k+1$ required to partition the $n$-dimensional Hamming space. Using a computer search, we improve upon the known upper bound for $n=9$ by showing that $13 \leq \chi_{\bar{2}}(9) \leq 14$.
\end{abstract}

\section{Introduction}
The \emph{$n$-hypercube} (or \emph{$n$-cube}), denoted by $Q_n$, is a graph whose vertex set corresponds to the set of all bit strings of length $n$. In the graph, two vertices are adjacent whenever their Hamming distance is one, i.e., the two bit strings differ in exactly one position. The \emph{$k$th power} of a graph $G$, denoted by $G^k$, is the graph on the same vertex set as $G$ with an edge between each pair of vertices at a distance at most $k$. In particular, $G^2$ is known as the \emph{square} of $G$. Originating from a study of scalability of optical networks, the problem of determining the chromatic number of the square of the $n$-hypercube has received considerable attention (see e.g.,~\cite{Wan1997, Ostergard2004, Das2006, Zhou2007, Fu2013}). In the context of coding theory, a proper vertex coloring of $Q_n^k$ corresponds to a partition of $\{0,1\}^n$ into binary codes of minimum distance at least $k+1$. We denote the chromatic number of $Q^k_n$ by $\chi_{\bar{k}}(n)$. Specifically, in what is to follow, we will focus on $\chi_{\bar{2}}(n)$.

Determining the exact value of $\chi_{\bar{k}}(n)$ has been deemed a difficult problem~\cite{Ngo2002}. For roughly 25 years, it was known that $13 \leq \chi_{\bar{2}}(8) \leq 14$, with explicit 14-colorings given independently by Hougardy~\cite{Ziegler2001} and Royle~\cite[Section~9.7]{Jensen2011}. Recently, Kokkala and Östergård~\cite{Kokkala2015} proved that $\chi_{\bar{2}}(8) = 13$ by constructing an explicit 13-coloring for $Q_8^2$. Given the result, the precise value of $\chi_{\bar{2}}(n)$ is known for all $n \leq 8$. For $n=9$, it is not difficult to derive that $ 13 \leq \chi_{\bar{2}}(9) \leq 16$. In fact, to the best of our knowledge, no tighter bounds are known. We improve upon the known upper bound by establishing that $13 \leq \chi_{\bar{2}}(9) \leq 14$. In particular, we demonstrate an explicit 14-coloring for the square of the 9-hypercube, found by a computer search.

\section{An explicit 14-coloring for \texorpdfstring{$Q_9^2$}{Q92}}
It is known that 16 colors suffice to properly color the vertices of $Q_9^2$ (for a gentle introduction, see~\cite{Rix2008}). Similarly, at least 13 colors are needed. In what is to follow, we present an explicit 14-coloring for $Q_9^2$, bringing the known upper bound from 16 down to 14. In other words, our construction establishes that $13 \leq \chi_{\bar{2}}(9) \leq 14$.

The claimed 14-coloring is given in Table~\ref{tbl:coloring}. Eight of the color classes are of size 40, two of size 35, two of size 33, one of size 31, and one of size 25. For convenience, we list the elements of each of the 14 color classes as integers from 0 to 511. We also provide a simple program\footnote{The source code is freely available at \url{https://bitbucket.org/laurij/sqhyp-check}} for verifying that the given coloring is indeed proper.

We will give details of the computational approach taken to discover the coloring in a full version of the manuscript.

\begin{table}[t]
\centering
\caption{A partition of the vertex set of $Q_9^2$ into color classes $C_1,\ldots,C_{14}$.}
\label{tbl:coloring}
\resizebox{\textwidth}{!}{
\begin{tabular}{|r|r|r|r|r|r|r|r|r|r|r|r|r|r|}
\hline
$C_1$ & $C_2$ & $C_3$ & $C_4$ & $C_5$ & $C_6$ & $C_7$ & $C_8$ & $C_9$ & $C_{10}$ & $C_{11}$ & $C_{12}$ & $C_{13}$ & $C_{14}$ \\
\hline
22 & 1 & 14 & 13 & 12 & 7 & 6 & 0 & 17 & 3 & 8 & 2 & 9 & 4 \\
25 & 24 & 27 & 23 & 19 & 10 & 11 & 15 & 40 & 16 & 20 & 5 & 21 & 18 \\
48 & 46 & 37 & 26 & 32 & 28 & 29 & 35 & 47 & 42 & 38 & 44 & 30 & 31 \\
63 & 55 & 50 & 34 & 58 & 36 & 33 & 53 & 66 & 45 & 78 & 59 & 39 & 43 \\
64 & 70 & 60 & 52 & 61 & 41 & 62 & 56 & 76 & 54 & 81 & 93 & 79 & 49 \\
71 & 85 & 77 & 57 & 75 & 51 & 72 & 73 & 95 & 69 & 101 & 99 & 96 & 74 \\
90 & 91 & 80 & 67 & 86 & 65 & 84 & 83 & 97 & 92 & 120 & 112 & 123 & 89 \\
107 & 105 & 87 & 68 & 88 & 82 & 98 & 94 & 116 & 115 & 127 & 126 & 155 & 118 \\
108 & 114 & 102 & 104 & 103 & 110 & 109 & 100 & 122 & 132 & 141 & 139 & 162 & 195 \\
117 & 124 & 121 & 111 & 113 & 125 & 119 & 106 & 131 & 143 & 151 & 158 & 177 & 253 \\
130 & 135 & 136 & 129 & 137 & 128 & 140 & 138 & 142 & 154 & 160 & 161 & 188 & 264 \\
133 & 146 & 145 & 134 & 144 & 149 & 147 & 148 & 152 & 181 & 185 & 183 & 197 & 277 \\
156 & 157 & 150 & 171 & 159 & 175 & 167 & 153 & 164 & 194 & 190 & 184 & 216 & 294 \\
169 & 168 & 163 & 172 & 165 & 186 & 170 & 166 & 178 & 215 & 196 & 200 & 233 & 301 \\
174 & 180 & 191 & 202 & 182 & 198 & 176 & 173 & 189 & 217 & 210 & 211 & 238 & 314 \\
179 & 187 & 203 & 208 & 206 & 205 & 193 & 199 & 201 & 235 & 234 & 212 & 247 & 335 \\
209 & 192 & 220 & 221 & 213 & 219 & 207 & 204 & 214 & 236 & 258 & 230 & 274 & 348 \\
223 & 222 & 224 & 229 & 226 & 227 & 218 & 225 & 231 & 240 & 287 & 237 & 312 & 352 \\
246 & 239 & 245 & 243 & 251 & 232 & 228 & 242 & 263 & 262 & 299 & 275 & 319 & 391 \\
248 & 241 & 250 & 254 & 252 & 244 & 249 & 255 & 265 & 289 & 300 & 276 & 326 & 409 \\
257 & 271 & 259 & 256 & 261 & 269 & 279 & 268 & 286 & 315 & 305 & 317 & 345 & 417 \\
270 & 278 & 260 & 267 & 266 & 272 & 280 & 273 & 304 & 316 & 327 & 321 & 355 & 436 \\
295 & 288 & 285 & 284 & 281 & 283 & 292 & 282 & 341 & 331 & 329 & 330 & 364 & 452 \\
296 & 313 & 298 & 293 & 291 & 290 & 303 & 297 & 344 & 338 & 346 & 356 & 396 & 466 \\
333 & 323 & 311 & 302 & 308 & 309 & 306 & 310 & 358 & 360 & 374 & 375 & 421 & 490 \\
339 & 332 & 328 & 307 & 320 & 318 & 325 & 322 & 365 & 373 & 385 & 384 & 427 & \\
340 & 336 & 350 & 337 & 351 & 324 & 334 & 349 & 371 & 392 & 398 & 397 & 438 & \\
354 & 357 & 353 & 342 & 361 & 343 & 347 & 359 & 404 & 403 & 408 & 430 & 448 & \\
377 & 362 & 367 & 378 & 366 & 363 & 369 & 368 & 411 & 413 & 423 & 434 & 467 & \\
382 & 383 & 372 & 381 & 370 & 376 & 380 & 379 & 426 & 418 & 469 & 479 & 478 & \\
395 & 388 & 399 & 402 & 390 & 387 & 386 & 389 & 439 & 449 & 493 & 491 & 509 & \\
400 & 394 & 410 & 405 & 412 & 406 & 393 & 415 & 463 & 462 & 496 & 497 & & \\
407 & 401 & 422 & 440 & 424 & 428 & 414 & 416 & 480 & 468 & 507 & 508 & & \\
420 & 419 & 425 & 447 & 431 & 441 & 437 & 435 & 505 & 487 & & & &\\
442 & 429 & 432 & 455 & 433 & 458 & 443 & 444 & 510 & 506 & & & &\\
445 & 446 & 450 & 460 & 451 & 465 & 464 & 459 &  &  &  &  &  & \\
454 & 457 & 453 & 475 & 461 & 476 & 477 & 470 &  &  &  &  &  & \\
456 & 471 & 473 & 482 & 474 & 485 & 483 & 472 &  &  &  &  &  & \\
481 & 486 & 492 & 489 & 484 & 498 & 488 & 494 &  &  &  &  &  & \\
495 & 504 & 499 & 500 & 503 & 511 & 502 & 501 &  &  &  &  &  & \\
\hline
\end{tabular}}
\end{table}

\bibliographystyle{abbrv}
\bibliography{sqhyp}

\end{document}